\documentclass[12pt]{article}
\usepackage{amsmath,amsfonts}
\usepackage{verbatim}
\usepackage{relsize}
\usepackage{enumerate}
\usepackage[inline]{enumitem}
\usepackage{color}
\normalbaselines
\newtheorem{teor}{Theorem}
\newtheorem{prop}[teor]{Proposition}
\newtheorem{lema}[teor]{Lemma}
\newtheorem{coro}[teor]{Corollary}
\newtheorem{rem}[teor]{Remark}

\newtheorem{ejem}[teor]{Example}

\newenvironment{demo}{\rm \trivlist \item[\hskip \labelsep{\it
      Proof}.]}{\nopagebreak \hfill $\square$ \endtrivlist}

\title{Rigidity results for complete spacelike submanifolds in plane fronted waves}

\author{Francisco J. Palomo${}^{a}$, Jos\'e A. S. Pelegr\'in${}^{b}$ and Alfonso
Romero${}^{c}$ \\[3mm]
${}^{a}$\normalsize Departamento de Matem\'atica Aplicada\\[-1mm] 
\normalsize Universidad de M\'alaga\\[-1mm]
\normalsize 29071 M\'alaga, Spain,\\[-1mm]
{\normalsize E-mail:} \small\texttt{fpalomo@uma.es}\\[2mm]
${}^{b}$\normalsize Departamento de Did\'actica de la Matem\'atica\\[-1mm]
\normalsize Universidad de Granada\\[-1mm]
\normalsize 18071 Granada, Spain,\\[-1mm]
{\normalsize E-mail:} \small\texttt{jpelegrin@ugr.es}\\[2mm]
${}^{c}$\normalsize Departamento de Geometr\'ia y Topolog\'ia\\[-1mm]
\normalsize Universidad de Granada\\[-1mm]
\normalsize 18071 Granada, Spain,\\[-1mm]
{\normalsize E-mail:} \small\texttt{aromero@ugr.es}\\}

\date{}

\begin{document}

\maketitle

\thispagestyle{empty}

\begin{abstract}
New rigidity results for complete non-compact 
spacelike submanifolds of arbitrary codimension in plane fronted waves are obtained. Under 
appropriate assumptions, we prove that a complete  
spacelike submanifold in these spacetimes 
is contained in a characteristic lightlike hypersurface. Moreover, 
for a complete codimension two extremal submanifold in 
a plane fronted wave we show sufficient conditions to guarantee 
that it is a (totally geodesic) wavefront.
\end{abstract}

\vspace*{5mm}


\section{Introduction}

In the search for exact solutions to Einstein's field equation it is usually assumed the existence of a certain symmetry. This symmetry is usually provided by a globally defined causal conformal vector field \cite{SK}. In particular, when this vector field is lightlike and parallel, the resulting solution is called a Brinkmann spacetime \cite{Br}. In this article, we will focus on a distinguished subfamily of Brinkmann spacetimes, namely, plane fronted waves \cite{FS}. A plane fronted wave is a Lorentzian manifold $(\overline{M}^{n+2},\langle \ , \ \rangle)$ where $\overline{M}^{n+2}=\mathbb{R}^{2}\times M^{n}$ with $(M^{n},g_{M})$ a (connected) Riemannian manifold and the Lorentzian metric 

\begin{equation}
\label{metrica}
\langle \ ,\ \rangle=  \mathcal{H}(u,x)du\otimes du +  du\otimes dv + dv \otimes du + g_{M},
\end{equation}

\noindent where $(u,v)$ are the natural coordinates on $\mathbb{R}^2$ and 
$\mathcal{H}(u,x)$ 
is a (independent of $v$) smooth function on $\mathbb{R}^{2}\times M^{n}$. 
The coordinate vector field $\partial_v := \partial / \partial v$ is lightlike and parallel. Endowing $\overline{M}^{n+2}$ with the time orientation defined by $\partial_v$, it 
becomes a spacetime. As $\overline{\nabla} u = \partial_v$, the coordinate 
$u : \overline{M}^{n+2} \longrightarrow \mathbb{R}$ plays the role of a 
quasi-time function \cite[Def. 13.4]{BEE} i.e., its gradient is everywhere 
causal and any causal segment $\gamma$ such that $u \circ \gamma$ constant 
is injective. In particular, the 
spacetime is causal \cite[p. 490]{BEE}. A plane fronted wave is 
foliated by the (characteristic) lightlike hypersurfaces $u = u_{0}$, $u_0 \in \mathbb{R}$. 

When $M^n = \mathbb{R}^2$ and $g_M$ is the usual Euclidean metric, the spacetime $\overline{M}^4$ is called a pp-wave (plane fronted wave with parallel propagating rays) \cite{BEE}. Plane fronted waves model (electromagnetic or gravitational) radiation propagating at the speed of light. Despite the fact that the study of gravitational waves goes back to Einstein and Rosen \cite{ER}, the experimental detection of gravitational waves \cite{Ab} has aroused widespread interest in plane fronted waves. 

In these ambient spacetimes we will study spacelike submanifolds of arbitrary codimension. Spacelike submanifolds of 
codimension greater than one became interesting from a physical viewpoint
when Penrose introduced the notion of trapped surface to study spacetimes' singularities \cite{Pe}. Namely, the existence of a trapped surface 
is a sign of the presence of a black hole. 

The original definition of trapped surface was given in terms of null expansions. This is related to the causal orientation of the mean curvature vector field, which allows the extension of the concept of trapped submanifold to arbitrary codimension \cite{Kr}. Indeed, let $\psi: \Sigma^k \longrightarrow \overline{M}^{n+2}$ be a 
spacelike submanifold of arbitrary codimension in an arbitrary 
spacetime. Denoting by $\overrightarrow{H}$ 
the mean curvature vector field of the submanifold and following the 
standard terminology in General Relativity (see \cite{Kr} and 
\cite{MS}), $\Sigma^k$ is said to be future trapped if $\overrightarrow{H}$ is timelike and 
future pointing everywhere (similarly for past trapped); weakly future trapped if $\overrightarrow{H}$ is causal and 
future pointing everywhere (similarly for weakly past trapped) and extremal if $\overrightarrow{H} = 0$.

In a plane fronted wave, the wavefronts $u = u_{0}$, $v = v_{0}$ are 
a distinguished family of codimension two extremal submanifolds. Therefore, 
it is natural to wonder under which conditions a complete codimension two extremal  submanifold is a wavefront. In \cite{CPR}, the authors answer this 
question for the compact case. Thus, in this paper we will focus on the 
non-compact case. Our main tools will be certain maximum principles as 
well as the parabolicity of the spacelike submanifold.
 
 Let us recall that a complete (non-compact) Riemannian 
manifold is parabolic if the only superharmonic functions
bounded from below that it admits are the constants 
(see, for instance, \cite{Ka}). From a 
physical perspective, parabolicity  
is equivalent to the recurrence of the Brownian motion on a Riemannian 
manifold \cite{Gr}. From a mathematical standpoint,  
complete spacelike surfaces with non-negative 
Gaussian curvature are parabolic \cite{H}. In arbitrary 
dimension there is no clear relation between parabolicity
and sectional curvature. Nevertheless, there exist sufficient conditions to 
ensure the parabolicity of a Riemannian manifold of arbitrary dimension 
based on its geodesic balls' volume growth \cite{AMR}. Along 
this article we will use that parabolicity is invariant 
under quasi-isometries \cite[Cor. 5.3]{Gr}.

This paper is organized as follows. Section \ref{sesu} introduces general notions 
on plane fronted waves and their spacelike submanifolds that will later 
allow us to proof our main results. Section \ref{semr} contains the main results of this paper,
we first state several conditions that guarantee that a non-compact 
spacelike submanifold is contained in a lightlike hypersurface $u = u_0$ 
(Theorems \ref{teoasy} and \ref{teopar}). From these theorems we can deduce 
certain non-existence results for spacelike hypersurfaces 
(Corollaries \ref{coronohyp} and \ref{corohyppar}) as 
well as for weakly trapped submanifolds 
(Theorem \ref{teoparwt}). Finally, 
we answer our initial question providing in 
Theorems \ref{teoparex}, \ref{teoyau} and \ref{teotcc} sufficient conditions 
for a complete non-compact codimension two extremal submanifold to be a wavefront. Indeed, Theorem \ref{teoparex} states that in a plane fronted wave the only codimension two parabolic extremal submanifolds 
with bounded (either from above or below) $u\circ \psi$ and $v\circ \psi$  
are the (necessarily parabolic) wavefronts. Furthemore, Theorem \ref{teoyau} asserts, under the same bound, that a complete codimension two extremal submanifold whose Ricci tensor 
verifies $\mathrm{Ric} \geq 0$
is also a wavefront.

\section{Set up}

\label{sesu}

\subsection{Plane fronted waves}

Let us consider a plane fronted wave $\overline{M}^{n+2}=\mathbb{R}^{2}\times M^{n}$ and let  $\mathcal{L}(M^{n})$ be the subspace of $\mathfrak{X}(\overline{M}^{n+2})$ consisting of the lifts to $\overline{M}^{n+2}$ of all vector fields on $M^{n}$. From now on, we will use the  the same symbol to denote a vector field in $\mathfrak{X}(M^{n})$ and its corresponding lift to $\mathcal{L}(M^{n})$. In a similar way, we simplify the notation using the same symbol for a function on $M^n$ and its corresponding lift to $\overline{M}^{n+2}$.
The Levi-Civita connection $\overline{\nabla}$ of a plane fronted wave was given in  \cite{CFS} as follows,

\begin{lema}\label{levicivita}
	For any $V, W\in \mathcal{L}(M^{n})\subset \mathfrak{X}(\overline{M}^{n+2})$, we have
	
	\begin{itemize}
		\item[\rm(i)] $\overline{\nabla}_{\partial_u}\partial_u=\frac{1}{2}\big( \partial_{u}\mathcal{H}\, \partial_{v}- \widetilde{\nabla} \mathcal{H}_u \big)= \frac{1}{2}\big(\overline{\nabla} \mathcal{H}- \widetilde{\nabla} \mathcal{H}_u\big)$, 
		\item[\rm(ii)] $\overline{\nabla}_{V} \partial_u=\overline{\nabla}_{\partial_u}V=\frac{1}{2}\, g_{M}(\widetilde{\nabla} \mathcal{H}_{u},V)\,\partial_v$,
		\item[\rm(iii)] $\overline{\nabla}_{V}W=\widetilde{\nabla}_{V}{W}$,		
		\item[\rm(iv)] $\overline{\nabla}_{\partial_v}\partial_v=\overline{\nabla}_{\partial_v}\partial_u=\overline{\nabla}_{\partial_u}\partial_v=\overline{\nabla}_{V} \partial_v=\overline{\nabla}_{\partial_v}V=0$,
	\end{itemize}
	where $\widetilde{\nabla}$ denotes the Levi-Civita connection of $(M^{n}, g_{M})$, $\widetilde{\nabla} \mathcal{H}_u$ is the gradient of $\mathcal{H}_u$ on $M^n$, $\mathcal{H}_{u}(x):=\mathcal{H}(u,x)$, and $\overline{\nabla} \mathcal{H}$ is the gradient of $\mathcal{H}$ on $\overline{M}^{n+2}$.  
\end{lema}

Every tangent vector $Z\in T_{(u,v,x)}\overline{M}^{n+2}$ admits the decomposition

\begin{equation}\label{161020A}
Z=\langle Z, \partial_{v}\rangle\,\partial_u +\big(\langle Z, \partial_{u}\rangle-\langle Z, \partial_{v}\rangle \mathcal{H}\big) \partial_v+d\pi_{M}(Z)
\end{equation}

\noindent where $\pi_{M}:\overline{M}^{n+2} \longrightarrow M^{n}$ is the natural projection. Thus, Lemma \ref{levicivita} implies that

\begin{equation}\label{120220A}
\overline{\nabla}_{Z}\partial_u=\frac{1}{2}\,\langle Z, \partial_v\rangle (\overline{\nabla} \mathcal{H} - \widetilde{\nabla}\mathcal{H}_{u})+ \frac{1}{2}\,g_{M}\big(\widetilde{\nabla}\mathcal{H}_{u}, d\pi_{M}(Z)\big)\,\partial_v.
\end{equation}

The non-necessarily vanishing components of the Riemann curvature tensor $\overline{\mathrm{R}}$ of the metric (\ref{metrica}) are
$$
\overline{\mathrm{R}}(V,\partial_u)\partial_u = -\frac{1}{2}\, \overline{\nabla}_V \overline{\nabla} \mathcal{H}_u  \quad \text{ and }\quad  \overline{\mathrm{R}}(V,\partial_u)W = \frac{1}{2}\, \widetilde{\mathrm{Hess}}(\mathcal{H}_u)(V,W)\,\partial_v,
$$
\noindent where $\widetilde{\mathrm{Hess}}$ stands for the Hessian tensor on $M^n$. Denoting by $\overline{\mathrm{Ric}}$ and $\mathrm{Ric}_M$ the Ricci tensors of $\overline{M}^{n+2}$ and $M^{n}$, respectively, we obtain that the non-necessarily vanishing components of $\overline{\mathrm{Ric}}$ are 
$$\mathrm{\overline{Ric}}(V,W)=\mathrm{Ric}_{M}(V,W) \ \ \text{and} \ \ \mathrm{\overline{Ric}}(\partial_u, \partial_u)=-\frac{1}{2}\widetilde{\Delta}\mathcal{H}_u,$$
\noindent where $\widetilde{\Delta}$ is the Laplacian on $M^{n}$. It can be easily seen that a plane fronted wave $(\overline{M}^{n+2}, \langle\ , \ \rangle)$ satisfies the timelike convergence condition (TCC), i.e., $\mathrm{\overline{Ric}}(T,T)\geq 0$ for every timelike vector $T$, if and only if 
\begin{equation}
\label{conditionsTCC}
\widetilde{\Delta}\mathcal{H}_u\leq 0 \ \ \text{and} \ \ \mathrm{Ric}_{M}\geq 0.
\end{equation}

\subsection{Spacelike submanifolds}

Let us now consider $\psi: \Sigma^k \longrightarrow \overline{M}^{n+2}$ a 
(connected) spacelike submanifold in plane fronted wave $(\overline{M}^{n+2}, \langle \ , \ \rangle)$. That is, $\psi$ is an immersion and the induced metric on $\Sigma^k$ is Riemannian. We will also denote the induced metric on $\Sigma^k$ by $\langle \ , \ \rangle$. For any vector field $V$ along the immersion $\psi$, we write $V^{\top}$ and $V^{\perp}$ for the tangent and normal parts along $\psi$, respectively. If we 
denote by $\mu:= u \circ \psi$ and $\nu := v \circ \psi$, the restrictions 
of the coordinate functions $u$ and $v$ to $\Sigma^{k}$, we obtain 
that their gradients are given by
\begin{equation}
\label{gru}
\nabla \mu = \partial_v^\top
\end{equation}
and
\begin{equation}
\label{grv}
\nabla \nu = \partial_u^\top - (\mathcal{H} \circ \psi) \partial_v^\top,
\end{equation}
Therefore, 
from these formulas we see how wavefronts are codimension two spacelike submanifolds 
such that the vector fields $\partial_u$ and $\partial_v$ are 
normal at every point.

Let us recall Gauss and Weingarten's formulas, given respectively by
\begin{equation}
\label{gauf}
\overline{\nabla}_X Y = \nabla_X Y + \mathrm{II}(X, Y),
\end{equation}
and
\begin{equation}
\label{weif}
\overline{\nabla}_X N = - A_N X + \nabla_X^\perp N,
\end{equation}
for 
$X, Y \in \mathfrak{X}(\Sigma^{k})$ and 
$N \in \mathfrak{X}^\perp (\Sigma^{k})$, where $\overline{\nabla}$ and $\nabla$ are the Levi-Civita connections of $\overline{M}^{n+2}$ and $\Sigma^k$, respectively, $\mathrm{II}$ 
is the second fundamental form and $A_N$ is the shape operator with respect to $N$ and $\nabla^\perp$ is the normal connection. The mean curvature vector field of the spacelike submanifold $\Sigma^{k}$ is

\begin{equation}
\label{mcv}
\overrightarrow{H} = \frac{1}{k} \ \mathrm{trace}_{\langle \ , \ \rangle}(\mathrm{II}).
\end{equation}

Taking into account that $\partial_{v}$ is a parallel vector field, we obtain using \eqref{gauf} and \eqref{weif} for 
$X \in \mathfrak{X}(\Sigma^{k})$

\begin{equation}
\label{wegaa}
0 = \overline{\nabla}_X \partial_v = \overline{\nabla}_X \partial_v^\top 
+ \overline{\nabla}_X \partial_v^\perp = \nabla_X \partial_v^\top 
+ \mathrm{II}(X, \partial_v^\top) - A_{\partial_v^\perp} X + \nabla_X^\perp \partial_v^\perp.
\end{equation}

Taking tangent and normal parts in (\ref{wegaa}) we get

\begin{equation}
\label{we1}
\nabla_X \partial_v^\top = A_{\partial_v^\perp} X
\end{equation}
and

\begin{equation}
\label{ga1}
\nabla_X^\perp \partial_v^\perp = - \mathrm{II}(X, \partial_v^\top).
\end{equation}

In particular, we have $\mathrm{div}(\partial_{v}^{\top})=\mathrm{trace}(A_{\partial_{v}^{\perp}})$.
Now, we can use (\ref{gru}), (\ref{grv}) and (\ref{we1}) to obtain 
that the Laplacians of $\mu$ and $\nu$ on $\Sigma^{k}$ are
\begin{equation}
\label{lau}
\Delta \mu = k \langle \overrightarrow{H}, \partial_v \rangle
\end{equation}
and
\begin{equation}
\label{lan}
\Delta \nu = \mathrm{div}(\partial_u^\top) 
- \partial_v^\top (\mathcal{H} \circ \psi) 
- k (\mathcal{H} \circ \psi) \langle \overrightarrow{H}, \partial_v \rangle.
\end{equation}

For a codimension two spacelike submanifold $\Sigma^{n}$ contained in 
a lightlike hypersurface $u = u_0$, the normal bundle is generated by the normal vector fields $\partial_{v}$ and $\eta=-\partial_{u}^{\perp}+\frac{1}{2}\langle \partial_{u}^{\perp}, \partial_{u}^{\perp}\rangle\partial_{v}$ which satisfy $\langle \eta, \eta \rangle=0$ and  $\langle \eta, \partial_{v} \rangle=-1$.
From (\ref{we1}) and \eqref{ga1}, we have 
\begin{equation}
\label{nabvlu}
A_{\partial_{v}}=0 \ \ \mathrm{and} \ \ \nabla^\perp \partial_v = 0.
\end{equation}
From \eqref{nabvlu}, the following formula holds for the second fundamental form
\begin{equation}
\label{sffco2}
\mathrm{II}(X,Y)=\langle A_{\partial_{u}^{\perp}}(X), Y\rangle \partial_{v},
\end{equation}
where $X,Y \in \mathfrak{X}(\Sigma^{n})$. Now, decomposing $\partial_u$ into its tangent and normal components and using (\ref{we1}) and (\ref{ga1}), we get from (\ref{120220A}) 
\begin{equation}
\label{derivco2}
\nabla_{W}\partial_u^\top =A_{\partial_u^{\perp}}(W) ,\ \nabla^{\perp}_W \partial_u^{\perp}=-\mathrm{II} (W, \partial_u^{\top})  +\frac{1}{2}\,\big\langle \nabla \mathcal{H}_{u}, (\pi_{M} \circ \psi)_{*}(W)\big\rangle\,\partial_v.
\end{equation}
Hence, from \eqref{lan}, \eqref{sffco2} and \eqref{derivco2}, the mean curvature vector of a spacelike submanifold $\Sigma^{n}$
contained in a lightlike hypersurface $u = u_0$ is given by (see \cite{CPR})

\begin{equation}
\label{hsuli}
\overrightarrow{H} =\frac{1}{n} \ \mathrm{trace}(A_{\partial_{u}^{\perp}}) \ \partial_v = \frac{1}{n} \ \mathrm{div}(\partial_u) \ \partial_v 
= \frac{1}{n} \ \Delta \nu \ \partial_v.
\end{equation}

The (codimension two) spacelike submanifolds given by $u = u_0$, $v = v_0$, for $u_0, v_0 \in \mathbb{R}$ are called wavefronts. From \eqref{nabvlu} and \eqref{derivco2} we 
have that each wavefront is totally geodesic. Each wavefront is clearly contained in a lightlike hypersurface $u = u_0$. As shown in \cite[Lemma 4.2]{CPR}, any codimension two spacelike submanifold of $\overline{M}^{n+2}$ contained in a lightlike hypersurface $u = u_0$ is locally isometric to $M^n$. Indeed, the following lemma extends this result.

\begin{lema}\label{lemaiso}
	Let $\psi:\Sigma^{n}\rightarrow \overline{M}^{n+2}$ be a codimension two spacelike submanifold in a plane fronted wave $\overline{M}^{n+2}$. If $\psi(\Sigma^n)$ is contained in a lightlike hypersurface $u = u_0$ or the lightlike vector field $2 \partial_u -(\mathcal{H}\circ \psi)\, \partial_v$ is normal to $\Sigma^n$ at any point, then, $\pi_{M}\circ \psi :\Sigma^{n} \to M^n$ is a local isometry.
\end{lema}

\begin{demo}
	From \eqref{161020A}, for any $Z \in T_p \Sigma^n$ we obtain
	$$ \langle d \psi_p (Z), d \psi_p (Z) \rangle = \langle Z , \partial_v \rangle \langle Z, 2 \partial_u - (\mathcal{H} \circ \psi)_p \partial_v \rangle + g_M (d (\pi_M \circ \psi)_p (Z), d (\pi_M \circ \psi)_p (Z)).$$
	
	\noindent Hence, the first addend on the right hand side vanishes under our assumptions.	
\end{demo}

\begin{rem}
\normalfont
 By means of equations (\ref{gru}) and (\ref{grv}) we obtain that 
 $2 \partial_u -(\mathcal{H}\circ \psi) \,\partial_v$ being normal to $\Sigma^n$  
 is equivalent to
$$
2\nabla \nu+ (\mathcal{H}\circ \psi )\, \nabla \mu=0.
$$
Therefore, if $2 \partial_u -(\mathcal{H}\circ \psi)\, \partial_v$ is normal to $\Sigma^n$ and 
$\psi(\Sigma^n)$ is contained in a lightlike hypersurface $u = u_0$, then 
$\psi(\Sigma^{n})$ is contained in a wavefront.

\end{rem}

\section{Main results}
\label{semr}

First,  let us obtain the next result for 
the restriction of the quasi time function $u$ to any spacelike 
submanifold.

\begin{prop}
\label{promax}
Let $\psi: \Sigma^k \longrightarrow \overline{M}^{n+2}$ be a 
spacelike submanifold in a plane fronted wave.

\begin{enumerate}[label=(\roman*)]

\item[\rm(i)] If $\langle \overrightarrow{H}, \partial_v \rangle > 0$, then 
the function $\mu$ attains no local maximum.

\item[\rm(ii)] If $\langle \overrightarrow{H}, \partial_v \rangle < 0$, then 
$\mu$ attains no local minimum.

\end{enumerate}

\end{prop}

\begin{demo}
Assume  $\mu$ attains a local maximum at  $p_{max}$. Therefore, at $p_{max}$ from \eqref{lau} we have
$\Delta \mu (p_{max}) = k \langle \overrightarrow{H}, \partial_v \rangle (p_{max}) 
\leq 0,$
\noindent contradicting the fact that 
$\langle \overrightarrow{H}, \partial_v \rangle > 0$. The 
proof of the second statement is analogous.
\end{demo}

As a direct consequence of Proposition \ref{promax}, we have the following result.

\begin{coro}\label{301121A}
Let $\psi: \Sigma^k \longrightarrow \overline{M}^{n+2}$ be a 
spacelike submanifold in a plane fronted wave with  $\langle \overrightarrow{H}, \partial_v \rangle > 0$ {\rm (}resp.  $\langle \overrightarrow{H}, \partial_v \rangle < 0${\rm)}. Then, there is no open subset $U\subset \Sigma^{k} $ such that $\psi(U)$ is contained in a lightlike hypersurface $u=u_{0}$.
\end{coro}

Another immediate consequence of Proposition \ref{promax} is 

\begin{coro}
\label{procpt}
Let $\psi: \Sigma^k \longrightarrow \overline{M}^{n+2}$ be a 
spacelike submanifold in a plane fronted wave. Then,

\begin{enumerate}[label=(\roman*)]

\item[\rm(i)] If $\mu$ attains a local maximum, then 
$\underset{\Sigma^k}{\mathrm{inf}} \langle \overrightarrow{H}, \partial_v \rangle 
\leq 0.$

\item[\rm(ii)] If $\mu$ attains a local minimum, then 
$\underset{\Sigma^k}{\mathrm{sup}} \langle \overrightarrow{H}, \partial_v \rangle 
\geq 0.$

\end{enumerate}

\end{coro}

Note that as an application of Proposition \ref{promax}
we reobtain the well known
non-existence of compact trapped submanifolds in a plane fronted wave 
(\cite{CPR}, \cite{MS}). 

Our next result ensures the constancy of the 
quasi time function's restriction to the spacelike submanifold under 
certain assumptions, equivalently, under which conditions the spacelike 
submanifold factorizes through the lightlike hypersurface $u = u_0$, for $u_0 \in \mathbb{R}$. 

\begin{teor}
	\label{teoasy}
	Let $\psi: \Sigma^k \longrightarrow \overline{M}^{n+2}$ be a 
	complete, orientable, non-compact spacelike 
	submanifold in a plane fronted wave.
	
	\begin{enumerate}[label=(\roman*)]
		
		\item[\rm(i)] If $\langle \overrightarrow{H}, \partial_v \rangle \geq 0,$ 
		$\mu \geq u_0,$ for $u_0 \in \mathbb{R}$ and 
		$$\lim_{r(p) \to \infty} \mu = u_0,$$ where $r(p) = d(p, o)$ is 
		the Riemannian distance function on $\Sigma^k$ from a fixed 
		point $o \in \Sigma^k$, then, $\mu = u_0$.
		
		\item[\rm(ii)] If $\langle \overrightarrow{H}, \partial_v \rangle \leq 0,$ 
		$\mu \leq u_0$ for $u_0 \in \mathbb{R}$ and 
		$$\lim_{r(p) \to \infty} \mu = u_0,$$ then, $\mu = u_0$.
	\end{enumerate} 
	
\end{teor}

\begin{demo}
	To prove the first statement, let us 
	define on $\Sigma^k$ the function $\widehat{\mu} := \mu - u_0$. From our 
	assumptions we have that $\widehat{\mu} \geq 0$ and 
	$$\lim_{r(p) \to \infty} \widehat{\mu} = 0.$$ Furthermore, using 
	\eqref{gru} we obtain 
	$$\langle \nabla \widehat{\mu}, \partial_v^\top \rangle = 
	|\partial_v^\top|^2 = |\nabla \mu|^2 \geq 0.$$
	Moreover, from our assumptions and \eqref{lau} we have
	$$\mathrm{div}(\partial_v^\top) = \Delta \mu = 
	k \langle \overrightarrow{H}, \partial_v \rangle \geq 0.$$
	
	Reasoning by contradiction, suppose that $\mu \neq u_0$. Hence, for 
	some $p \in \Sigma^k$ we have $\widehat{\mu}(p) > 0$, so that 
	$\widehat{\mu} \neq 0$. Thus, we can apply 
	\cite[Thm. 2.2]{ACN} to conclude that
	
	$$\langle \nabla \widehat{\mu}, \partial_v^\top \rangle = 
	|\nabla \mu|^2 = 0.$$
	
	Therefore, $\mu$ must be constant on $\Sigma^k$ and since $\mu$ is 
	asymptotic to $u_0$ we obtain $\mu = u_0,$ reaching a contradiction. The 
	proof of (ii) is analogous defining on $\Sigma^k$ the function 
	$\widetilde{\mu}:= u_0 - \mu$ and using the vector field 
	$-\partial_v^\top$ instead of $\partial_v^\top$ in the computations.
\end{demo}

In particular, for a spacelike submanifold 
$\psi: \Sigma^k \longrightarrow \mathbb{L}^{n+2}$ in the Lorentz-Minkowski spacetime 
$\mathbb{L}^{n+2}$, calling Beltrami's equation $\Delta \psi = k \overrightarrow{H}$, 
we derive from Theorem \ref{teoasy} the next result.

\begin{coro}
	\label{coromink}
	
Let $\psi: \Sigma^k \longrightarrow \mathbb{L}^{n+2}$ be a 
	complete, orientable, non-compact spacelike 
	submanifold in $\mathbb{L}^{n+2}$ with $\psi=(\psi_{0}, \psi_{1}, \dots , \psi_{n+1})$. Assume that there is $j\in \{1,\dots , n+1\}$ such that
	$\Delta (\psi_{j} - \psi_{0}) \geq 0$ and $\psi_{j}- \psi_{0}\geq u_{0}$, for $u_0 \in \mathbb{R}$. If
		$$\lim_{r(p) \to \infty} (\psi_{j}-\psi_{0}) = u_0,$$ where $r(p) = d(p, o)$ is 
		the Riemannian distance function on $\Sigma^k$ from a fixed 
		point $o \in \Sigma^k$, then $\psi$  factorizes through the lightlike hyperplane $x_{j}- x_{0}= u_{0}$.
\end{coro}

\begin{rem}
	\normalfont
	For the case of a spacelike hypersurface ($k = n+1$), the orientability assumption 
	in Corollary \ref{coromink} can be dropped. Indeed, the orientation and time orientation of $\mathbb{L}^{n+2}$ guarantees the orientability of $\Sigma^{n+1}$. In general, if the ambient spacetime $\overline{M}^{n+2}$ is orientable, then every spacelike hypersurface 
	in $\overline{M}^{n+2}$ is also orientable.

\end{rem}

\begin{rem}
	
	\normalfont
	
	Note that for the codimension two case, under the assumptions of Lemma \ref{lemaiso},  $\Sigma^n$ is locally isometric 
	to $M^n$. Therefore, if the plane fronted wave verifies the TCC, then 
	the Ricci curvature of $\Sigma^n$ will be non-negative \eqref{conditionsTCC}. In the particular 
	case where $n=2$, this fact means that the Gaussian curvature of $\Sigma^2$ is 
	non-negative. Hence, if $\Sigma^2$ is also complete, it must be parabolic by 
	\cite{H}. More generally, in a plane fronted wave 
	$\overline{M}^{n+2} = \mathbb{R}^2 \times M^n$ with $M^n$ parabolic, 
	a simply connected complete 
	codimension two spacelike submanifold with $\mu = u_0$ or everywhere orthogonal to $2 \partial_u - \mathcal{H} \partial_v$ is also 
	parabolic due to Lemma \ref{lemaiso} 
	and the fact that parabolicity is 
	invariant under	quasi-isometries \cite[Cor. 5.3]{Gr}.
\end{rem}

Furthermore, for parabolic spacelike submanifolds we have the next theorem.

\begin{teor}
\label{teopar}
Let $\psi: \Sigma^k \longrightarrow \overline{M}^{n+2}$ be a 
parabolic
spacelike submanifold in a plane fronted wave. If 
$\langle \overrightarrow{H}, \partial_v \rangle \geq 0$ 
{\rm(}resp., $\langle \overrightarrow{H}, \partial_v \rangle \leq 0${\rm)}  and 
$\sup_{\Sigma^k}  \ \mu < +\infty$ 
{\rm(}resp., $\inf_{\Sigma^k} \ \mu > -\infty${\rm)}, then $\mu = u_0$. In addition, if 
$\Sigma^n$ is a simply-connected codimension two submanifold, then the wavefront 
is parabolic.
\end{teor}

\begin{demo}
From \eqref{lau}, if $\langle \overrightarrow{H}, \partial_v \rangle \geq 0$ holds, 
$\mu$ would be a superharmonic function on a parabolic manifold. Thus, if it is 
bounded from below, it must be constant. We can prove the other case in a similar way.

To conclude, the last statement is a consequence of the fact that in the 
codimension two case, if $\mu$ is constant and $\Sigma^n$ is simply connected, then 
it is globally isometric to the wavefront \cite[Prop. 4.6]{CPR}.
\end{demo}

As an immediate consequence of Theorem \ref{teopar} we deduce that  
plane fronted waves do not admit compact spacelike hypersurfaces 
which are everywhere non-contracting or 
non-expanding.

\begin{coro}
\label{coronohyp}
In a plane fronted wave there are no compact spacelike hypersurfaces with signed 
mean curvature.
\end{coro}

Moreover, if we consider a spacelike hypersurface in a plane fronted wave 
and choose its normal unitary vector $N$  such that $\langle N, \partial_v \rangle < 0$, we obtain 
from Theorem \ref{teoasy} the next result for parabolic spacelike 
hypersurfaces (compare with \cite[Thm. 4]{PRR} and \cite[Thm. 4]{VL}).

\begin{coro}
\label{corohyppar}
In a plane fronted wave there are no parabolic spacelike hypersurfaces with 
non-positive mean curvature and $\mu$ bounded from above 
nor parabolic spacelike hypersurfaces with 
non-negative mean curvature and $\mu$ bounded from below.
\end{coro}

Also, by means of Theorem \ref{teopar} we can deduce the 
following non-existence result for 
parabolic weakly trapped submanifolds in a plane fronted wave.

\begin{teor}
\label{teoparwt}
In a plane fronted wave there are no codimension two parabolic weakly future trapped 
{\rm(}resp., weakly past trapped{\rm)} submanifolds 
with $\inf_{\Sigma^n} \ \mu > -\infty$ and 
$\sup_{\Sigma^n}  \ \nu < +\infty$ 
{\rm(}resp., $\sup_{\Sigma^n}  \ \mu < +\infty$  and 
$\inf_{\Sigma^n} \ \nu > -\infty${\rm)}.
\end{teor}

\begin{demo}
Let us assume the existence of a codimension two 
parabolic weakly future trapped submanifold in a 
plane fronted wave with $\mu$ bounded from below and $\nu$ bounded from above. 
From Theorem \ref{teopar} we obtain that it satisfies $\mu = u_0$. Moreover, 
from \eqref{hsuli} we obtain $\Delta \nu > 0$, which combined with the upper 
bound of $\nu$ and the parabolicity of the submanifold implies that $\nu$ 
is constant and therefore, the submanifold would be extremal, 
reaching a contradiction. We can prove the other case in an analogous manner.
\end{demo}

Using these ideas we get the next rigidity result for the 
extremal case in codimension two, which extends \cite[Thm. 4.10]{CPR} 
to the non-compact case.

\begin{teor}
\label{teoparex}
In a plane fronted wave the only codimension two parabolic extremal submanifolds 
with $\mu$ and $\nu$ bounded {\rm(}either from above or below{\rm)}
are the {\rm(}necessarily parabolic{\rm)} wavefronts. 
\end{teor}

\begin{demo}
From \eqref{lau} we obtain that $\mu$ is a bounded harmonic function in a 
parabolic Riemannian manifold. Thus, $\mu = u_0$ on $\Sigma^n$. In addition, 
from \eqref{hsuli} we have that since $\nu$ is also a bounded harmonic function 
$\nu = v_0$.
\end{demo}

In addition, using Yau's result obtained in \cite[Cor. 1]{Ya} we 
can extend this result to 
codimension two spacelike submanifolds with non-negative Ricci curvature.

\begin{teor}
\label{teoyau}
In a plane fronted wave the only complete codimension two extremal submanifolds 
with $\mu$ and $\nu$ bounded {\rm(}either from above or below{\rm)} and whose Ricci tensor 
verifies $\mathrm{Ric} \geq 0$
are the wavefronts.
\end{teor}

To conclude, we can provide another extension to the non-parabolic 
case by means of Theorem \ref{teoasy} as follows.

\begin{teor}
\label{teotcc}
Let $\psi: \Sigma^n \longrightarrow \overline{M}^{n+2}$ be a 
complete, orientable, non-compact extremal 
submanifold in a plane fronted wave that satisfies the {\rm TCC}. If $\mu$ and $\nu$ are bounded 
{\rm(}either from above or below{\rm)} by $u_0, v_0 \in \mathbb{R}$, respectively, 
and $$\lim_{r(p) \to \infty} \mu = u_0,$$ then $\Sigma^n$ is a wavefront.
\end{teor}

\begin{demo}
From Theorem \ref{teoasy}, we obtain that $\mu = u_0$. Moreover, 
if $\overline{M}^{n+2}$ satisfies the TCC, \cite[Prop.12]{CPR} guarantees that 
$\Sigma^n$ has non-negative Ricci curvature. Now, using \eqref{hsuli} and our 
assumptions we obtain that $\nu$ is a bounded harmonic function on a Riemannian 
manifold with non-negative Ricci curvature. Thus, it is constant by 
\cite[Cor. 1]{Ya}.
\end{demo}

\begin{ejem}
	
	\normalfont

A remarkable family of codimension two spacelike submanifolds in a plane fronted wave $\overline{M}^{n+2}$ satisfying $\mu = u_{0}$ are the spacelike graphs introduced in 
\cite[Sec. 4]{CPR}. Namely, given a smooth function 
$h\in C^{\infty}(\Omega)$, where $\Omega$ is an open domain in 
$M^n$ and $u_0\in\mathbb{R}$, the graph

\begin{equation}
\label{150220A}
\Sigma^n_{u_0}(h)=\big\{\,\big(u_0 , h(x),x\big)\in \mathbb{R}\times\mathbb{R}\times M^n \; :\; x\in \Omega\,\big\},
\end{equation}

\noindent defines a codimension two spacelike submanifold, for any $h$ and $u_0$.

Note that, as a direct consequence of \eqref{hsuli}, the mean curvature of the graph 
$\Sigma^n_{u_0}(h)$ is $\overrightarrow{H}=(\Delta h /n)\,\partial_v$.
Clearly, each of these graphs is a 
wavefront if and only if $h$ is constant. Thus, the assumptions 
in Theorems \ref{teoparex}, \ref{teoyau} and \ref{teotcc} cannot be weakened, since we 
would easily find counterexamples within this family of spacelike graphs.
\end{ejem}

\section*{Acknowledgements} The  first author is partially
supported by Spanish MICINN project\linebreak PID2020-118452GB-I00. The last two authors  by Spanish MICINN project PID2020-116126GB-I00. The first and the third authors are also supported by Andalusian and ERDF project P20$_{-}$01391.
Research partially supported by the ``Mar\'{\i}a de Maeztu'' Excellence
Unit IMAG, reference CEX2020-001105-M, funded by
MCIN-AEI-10.13039-501100011033.


\begin{thebibliography}{99}

\bibitem{Ab} B.P. Abbott et al., Observation of gravitational waves from a binary black hole merger, \emph{Phys. Rev. Lett.}, \textbf{116} (2016), 1-–16.

\bibitem{ACN} L.J. Al\'ias, A. Caminha and Y. do Nascimento, A maximum 
principle at infinity with applications to geometric vector 
fields, \emph{J. Math. Anal. Appl.}, \textbf{474} (2019), 242--247.

\bibitem{AMR} L.J. Al\'ias, P. Mastrolia and M. Rigoli, 
\emph{Maximum principles and geometric applications}, Springer, 2016.

\bibitem{BEE} J.K. Beem, P.E. Ehrlich and K.L. Easley, 
\emph{Global Lorentzian geometry},
Monographs Textbooks Pure Appl. Math. 202, Dekker Inc., New York, 1996.

\bibitem{Br} H. Brinkmann, Einstein spaces which are mapped 
conformally on each other, {\emph Math. Ann.}, \textbf{94} (1925), 119--145.

\bibitem{CFS} A. Candela, J.L. Flores and M. S\'anchez, 
On general plane fronted waves. Geodesics, 
\emph{Gen. Relat. Gravit.}, \textbf{4} (2003), 631--649.

\bibitem{CPR} V.L. C\'anovas, F.J. Palomo and A. Romero, Mean curvature 
of spacelike submanifolds in a Brinkmann spacetime, 
\emph{Classical Quant. Grav.}, \textbf{38} (2021), 1--18.

\bibitem{ER} A. Einstein and N. Rosen, On gravitational waves, \emph{J. Franklin Inst.}, \textbf{223} (1937), 43--54.

\bibitem{FS} J.L. Flores and M. S\'{a}nchez, \textit{On the 
geometry of pp-wave type spacetimes},  Analytical and 
numerical approaches to mathematical relativity, Lecture 
Notes in Phys., {\bf 692}, Springer, Berlin, (2006), 79--98.

\bibitem{Gr} A. Grigor'yan, Analytic and geometric background of recurrence 
and non-explosion of the Brownian motion on Riemannian manifolds, 
\emph{B. Am. Math. Soc.}, \textbf{36} (1999), 135--249.

\bibitem{H} A. Huber, On subharmonic functions and differential geometry in the large, 
\emph{Comment. Math. Helv.}, \textbf{32} (1958), 13--72.

\bibitem{Ka} J.L. Kazdan, Parabolicity and the Liouville property on complete 
Riemannian manifolds, \emph{Aspects of Math.}, \textbf{10} (1987), 153--166.

\bibitem{Kr} M. Kriele, \emph{Spacetime: Foundations of General Relativity 
and Differential Geometry}, Springer Science and Business Media, 1999.

\bibitem{MS} M. Mars and J.M.M. Senovilla, Trapped surfaces and 
symmetries, \emph{Classical Quant. Grav.}, \textbf{20} (2003), L293--L300.

\bibitem{PRR} J.A.S. Pelegr\'in, A. Romero and R.M. Rubio, 
On maximal hypersurfaces in Lorentz manifolds
admitting a parallel lightlike vector field, \emph{Classical Quant. Grav.}, 
\textbf{33} (2016), 055003(1--8).

\bibitem{Pe} R. Penrose, Gravitational collapse and space-time singularities, 
\emph{Phys. Rev. Lett.}, \textbf{14} (1965), 57.

\bibitem{SK} H. Stephani, D. Kramer, M. MacCallum, C. Hoenselaersand and E. Herlt, \emph{Exact Solutions of Einstein’s Field Equations}, Cambridge University Press, 2003.

\bibitem{VL} M.A.L. Vel\'asquez and H.F. de Lima, Complete spacelike 
hypersurfaces immersed in pp-wave
spacetimes, \emph{Gen. Relat. Gravit.}, \textbf{52} (2020), 41 (1--18).

\bibitem{Ya} S.T. Yau, Harmonic functions on complete Riemannian manifolds, 
\emph{Comm. Pure and Applied Math.}, \textbf{28} (1975), 201--228.

\end{thebibliography}
\end{document}